\documentclass[12pt]{article}
\usepackage{amsmath, amsthm, euscript}
\newtheorem{theorem}{Theorem}[section]

\newtheorem{lemma}{Lemma}[section]

\newtheorem{remark}{Remark}[section]

\usepackage{graphicx, amsmath}

\begin{document}
\normalsize
\title{A proof of a conjecture in the Cram\'{e}r-Lundberg model with investments}

\maketitle

Shimao Fan, Sheng Xiong, Wei-Shih Yang
\vskip 12pt
Department of Mathematics

        Temple University

        Philadelphia, PA 19122-6094

\vskip 12pt
Email: shimao.fan@temple.edu, sheng@temple.edu, yang@temple.edu

\vskip 12pt


\vskip12pt
AMS classification Primary: 60J70, 91B30

\begin{abstract}
In this paper, we discuss the  Cram\'{e}r-Lundberg model with investments, where the price of the invested risk asset follows a geometric Brownian motion with drift $a$ and volatility $\sigma> 0.$ By assuming there is a cap on the claim sizes, we prove that the  probability of ruin has at least an algebraic decay rate if $2a/\sigma^2 > 1$. More importantly, without this assumption, we show that the  probability of ruin is certain for all initial capital $u$, if $2a/\sigma^2 \le 1$.
\end{abstract}

\section{Introduction}

\setcounter{equation}{0}
 In the classical Cram\'{e}r-Lundberg model, if the claim sizes have finite exponential moments, then it is well-known that the ruin probability decays exponentially as the initial surplus increases; see for instance the books by Asmussen \cite{Asm} and Embrechts et. al \cite{EKM}. For the case of
heavy-tailed claims there also exists numerous results in the literature. However, when the insurance company invests in a risky asset, for example a stock, whose price is described by a geometric Brownian motion with drift $a>0$ and volatility $\sigma> 0,$ then the probability of ruin either decays algebraically  as the initial surplus increases or the ruin is certain, provided the claim size is exponentially distributed. This result was shown by Frolova et. al \cite{FKP}. Under the assumption that the claim size distributions have moment generating functions defined on a neighborhood of the origin, Constantinescu and Thommann \cite{corina} proved that if the probability of ruin decays as the initial capital $u\rightarrow \infty,$ then $\rho=\frac{2a}{\sigma^2}>1$, and that if $1<\rho < 2$, then the probability of ruin decays algebraically as the initial capital  $u\rightarrow \infty$. Furthermore, they conjectured that if $\rho \le 1$, then the ruin probability $\psi(u) = 1$ for all $u\ge0$.

In this paper, our main goal is to prove that the conjecture is true. This work was motivated by a paradox of risk without the possibility of reward discussed by Steele \cite{Steele}. In the setting of this paradox of risk, the price of a risky asset is modeled by a geometric Brownian motion with an expected return rate $a$.  Steele pointed out that if $\rho<1$, the price of the risky asset approaches to zero with probability one, despite the fact that the expected value goes to positive infinity at an exponential rate. We observe that if the price of our risky asset is very close to zero, then even a small jump will trigger the ruin. Similarly, if the price of the risky asset drops below a threshold with probability one and if there is a positive probability that the price of the risky asset may have jumps larger than the threshold, then the ruin occurs almost surely. If the jump is modeled by a compound Poisson process, then this leads to the conjecture that is discussed in this paper.

We first recall the Cram\'{e}r-Lundberg model with investments. The risk process is given by
\begin{eqnarray}\label{C_L}
X_t=X_0+\int_0^t a X_s dt + \int_0^t \sigma X_sdW_s +ct-\sum_{j=1}^{N(t)} \xi_j ,
\end{eqnarray}
or
\begin{eqnarray}\label{C_L1}
 dX_t=(aX_t+c)dt+\sigma X_t dW_t-dP_t,
\end{eqnarray}
where $W_t$ is the Wiener process (standard Brownian motion), $N(t)$ is a Poisson process with parameter $\lambda$, and the claim sizes
 $\xi_i$;         $i=1, 2, 3,...$, are independent, identically distributed random variables, having the density function $p(x)$, with positive mean $\mu$ and finite variance. $c$ is the fixed rate of premium and $X_0$ is the initial capital. $P_t=\sum_{j=1}^{N(t)} \xi_j.$ The capital $X_t$ is continuously invested in a risky asset, with relative price increments $dX_t = aX_tdt+\sigma X_t dW_t,$ where $ a>0$ and $\sigma>0$ are the drift and volatility of the returns of the asset.\smallskip

Our paper is organized as follows. By assuming there is a cap on the  claim size, in Section 2, we prove two important results that (1) the  probability of ruin has at least an algebraic decay rate if $2a/\sigma^2 > 1$ and  (2) the price of the risky asset will drop below a threshold with probability one for all initial capital $X_0=u$, if $2a/\sigma^2 \le 1$. In Section 3, we prove that the conjecture is true by coupling the stochastic processes with and without the assumption on the claim sizes.



\section{Ruin Probability With A Cap On the Claim Size}
 We will assume the claim size is bounded by a constant $M>0$ through the entire section. In insurance, $M$ can be understood as the limit or cap of a policy.
\noindent Let
$T_{u^*}=\inf \{t > 0;\; X_t < u^* \}$
be the first time that $X_t < u^*$, and let
$$\psi_{u^*}(u)=P( T_{u^*} <\infty \;| X_0=u)$$
be the probability of ruin, where $0 \le u^*<u$. If $u^*=0$, we denote the probability of ruin by $\psi(u)$. We will discuss the  probability of ruin on the Cram\'{e}r-Lundberg model with investments based on (1) $\rho>1,$ (2) $\rho=1$ and (3) $\rho<1$.
\noindent We first prove the following
\begin{lemma}\label{Decreasing}
Let $X_t$ be a stochastic process that satisfies (\ref{C_L1}), if $0\le v\le u.$ then $$\psi(v)\ge \psi(u).$$
\end{lemma}
 Proof. We first derive a strong solution for (\ref{C_L1}). Let $Y_t=\exp\{(\frac{\sigma^2}{2}-a)t-\sigma W_t\}$. By It\^{o}'s formula \cite{IW}, $dX_tY_t=X_tdY_t+Y_tdX_t+dX_tdY_t$, and simple calculation yields $dX_tY_t=d{V_t}^u,$ where ${V_t}^u=u+c\int_0^tY_s\;ds-\int_0^tY_s\;dP_s.$
 Integrating both sides, we have $X_tY_t={V_t}^u$. Hence
\begin{eqnarray}
X_t=Y_t^{-1}{V_t}^u \label{eq:strong solution}
\end{eqnarray} is a strong solution of (\ref{C_L}) and (\ref{C_L1}) with initial condition  $X_0=u.$

 Next we define $Z_t=Y_t^{-1}{V_t}^v,$ then $Z_t \le X_t,\; \forall \;t\ge 0,$ since $0\le v\le u.$ Hence
\begin{eqnarray}
\psi(u)=P(X_t<0, \mbox{ for some } \; 0< t < \infty\;| X_0=u) \nonumber \\
\le P(Z_t<0, \mbox{ for some } \; 0< t < \infty\;| Z_0=v). \nonumber
\end{eqnarray}
Note that $Z_t$ also satisfies (\ref{C_L1}) with initial condition  $Z_0=v.$
 Hence
$$P(Z_t<0, \mbox{ for some } \; 0< t < \infty\;|Z_0=v)=\psi(v).$$
 Therefore
$$\psi(v)\ge \psi(u).$$

Our main tool is It\^{o}'s formula for semimartingales with a jump part.  Let $t_1 < t_2 < t_3 < ...$ be the times where the Poisson process $N(t)$ has a jump discontinuity.  Then the jump discontinuities for $P_t$ are also at $t_i$ with jump size $\xi_i$.  Following the notations on P. 43 \cite{IW}, for $t >0$, and a Borel set $U$ in $R$, we let
$$ N_p( (0, t] \times U)= \sharp \{i; t_i \le t,  \xi_i \in U \}. $$
Then $ N_p((0, t] \times U)$ defines a random measure $ N_p(dtdx)$ on the Borel $\sigma$-algebra on $[0, \infty) \times R$.  Note that
\begin{eqnarray}
N_p(dtdx)=\sum_{i=1}^{\infty}\delta_{t_i}(dt)\delta_{\xi_i}(dx), \label{eq:dirac}
\end{eqnarray}
where $\delta_{t_i}$ is the Dirac $\delta$-function centered at $t_i$ (probability measure concentrated at one point $t_i$).
It follows that
\begin{eqnarray}
\int_{0}^{t}\int_{0}^{\infty}f(s, x)N_p(dsdx)=\sum_{i; t_i\le t}f(t_i, \xi_i),  \label{eq:dirac sum}
\end{eqnarray}
and therefore
\begin{eqnarray}
\int_{0}^{t}\int_{0}^{\infty}xN_p(dsdx)=\sum_{i; t_i\le t}\xi_i=P_t.  \label{eq:dirac sum}
\end{eqnarray}

  It is well-known, see e.g. P. 60 and P. 65 \cite{IW}, that there exists a continuous process $\hat {N}_p((0, t] \times U)$ such that
  \begin{eqnarray}
\tilde {N}_p((0, t] \times U)=N_p((0, t] \times U)-\hat {N}_p((0, t] \times U), \label{eq:poisson mart}
\end{eqnarray}
is a martingale. In our case $$\hat {N}_p((0, t] \times U)=E[N_p((0, t] \times U)]. $$  $E[N_p((0, t] \times U)]$ defines a measure, $n_p(dtdx)$, called the mean (intensity) measure of $ N_p(dtdx)$ and it is given by $n_p(dtdx)=\lambda p(x) dtdx $.

The equation (\ref{C_L}) can be rewritten as
\begin{eqnarray}
X_t=X_0+\int_0^t a X_s dt + \int_0^t \sigma X_sdW_s +ct-\int_{0}^{t}\int_{0}^{\infty}xN_p(dsdx). \label{eq: semimartingale equation}
\end{eqnarray}
By (\ref{eq:strong solution}), the equation (\ref{eq: semimartingale equation}) has a strong solution for each fixed initial condition and it is a semimartingale by Definition 4.1, P. 64 \cite{IW}.

By (\ref{eq:strong solution}) and direct calculation, we have
\begin{eqnarray}
X_{t+s}=\bar{Y}_t^{-1}X_s+\bar{Y}_t^{-1}\int_{0}^{t}c\bar{Y}_u du - \bar{Y}_t^{-1}\int_{0}^{t}\bar{Y}_u d\bar{P}_u, \label{eq: Markov equation}
\end{eqnarray}
where
\begin{eqnarray}
\bar{Y}_t&=&e^{-(a-\frac{\sigma^2}{2})t- \sigma \bar{W}_t},\\
\bar{W}_t&=&W_{t+s}-W_{s},\\
\bar{P}_t&=&P_{t+s}-P_{s}.
\end{eqnarray}

Note that
$\bar{W}_t$ and $\bar{P}_t$ are independent of $\{X_v; 0\le v \le s \}$ and therefore
given $\{X_v; 0\le v \le s \}$, $X_{t+s}$ depends on $X_s$ only.  This implies that $X_t$ is a Markov process.  Moreover, since $ \bar{W}_t=W_{t+s}-W_{s}$ and $W_t$ have the same distribution, and $\bar{P}_t=P_{t+s}-P_{s}$ and $P_t$ have the same distribution, we have
\begin{eqnarray}
P(X_{t+s} \in U | X_s=x)= P(X_{t} \in U | X_0=x),
\end{eqnarray}
for all $t > 0$, and all Borel sets $U$.  Therefore, $X_t, t \ge 0$ is a Markov process with a stationary transition function.  Since the sample paths of $X_t$ are right continuous with left limits, $X_t, t \ge 0$ is a strong Markov process.

\begin{theorem}
Consider the model given by (\ref{C_L}) and assume that  $\sigma> 0$,\;$\rho>1$ and $c> \lambda\mu$. Then $$\psi(u)\le\left(\frac{M}{u}\right)^{\rho-1} \;\; \forall \;u\ge M.$$
\end{theorem}

\begin{remark} This theorem shows that the probability of ruin has at least an algebraic decay rate if $2a/\sigma^2 > 1$. In fact, we obtain a slightly stronger result in the proof below:$$ \psi_M(u)\le\left(\frac{M}{u}\right)^{\rho-1} \;\; \forall \;u\ge M. $$

\end{remark}
Proof. Let $F(x)=x^{1-\rho}, x>0.$ Applying It\^{o}'s formula \cite{IW}, we have
\begin{align*}
F(X_t)-F(X_0)&=\int_0^t(1-\rho)(X_s)^{-\rho}(aX_s+c)\;ds+
\int_0^t(1-\rho)(X_s)^{-\rho}\sigma X_s dW_s\\
&+\frac{1}{2}\int_0^t(1-\rho)(-\rho)(X_s)^{-\rho-1}\sigma^2{X_s}^2\;ds\\
&+\int_0^{t^+}\int_0^{M}(X_{s^-}-x)^{1-\rho}-(X_{s^-})^{1-\rho}\;N_p(dsdx),
\end{align*}
Hence
\begin{eqnarray}
F(X_t)&=&F(X_0)+\int_0^t(1-\rho)(X_s)^{-\rho}(aX_s+c)\;ds+
\mbox{ mart. } \nonumber\\
&+&\frac{1}{2}\int_0^t(1-\rho)(-\rho)(X_s)^{-\rho-1}\sigma^2{X_s}^2\;ds \nonumber \\
&+&\int_0^{t^+}\int_0^{M}(X_{s^-}-x)^{1-\rho}-(X_{s^-})^{1-\rho}\; \tilde{N}_p(dsdx) \nonumber \\
 &\le& F(X_0)+ \mbox{ mart. }+c(1-\rho)\int_0^t (X_s)^{-\rho} ds \nonumber \\
&+&\int_0^{t^+}\int_0^{M}(1-\rho)(X_{s^-})^{-\rho}(-x)\lambda p(x)dxds, \label{eq:mart ineq rho>1} \\
&=&F(X_0)+ \mbox{ mart. }+(1-\rho)(c-\lambda\mu)\int_0^{t} (X_{s})^{-\rho}ds,
\end{eqnarray}
here, and through-out  this paper, $\mbox{ mart. }$ denotes a martingale at time $t$.
The inequality (\ref{eq:mart ineq rho>1}) holds because $$(X_{s^-}-x)^{1-\rho}-(X_{s^-})^{1-\rho}\le (1-\rho)(X_{s^-})^{-\rho}(-x),\; \forall X_{s^-}\ge M.$$

 Now we consider the process $X_t$ on $[M, n)$, where $n$ is an integer ($>M$). Let $$\tau_n=\inf\{t>0:\; X_t \not \in [M, n)\}$$ be the first exit time from the interval $ [M, n)$. By the Optional Stopping Theorem, it follows that
\begin{eqnarray}\label{Expect}
E[F(X_{\tau_n})]\le E[F(X_0)].
\end {eqnarray}
Since $\xi_j>0$ for all $j=1, 2,\ldots$, we have $X_{\tau_n}=n$ or $X_{\tau_n}< M$. Moreover, since $F(x)$ is decreasing, we have
$$E[F(X_{\tau_n})]\ge \frac{1}{M^{\rho-1}}P(X_{\tau_n}< M\;|X_0=u) + \frac{1}{n^{\rho-1}}P(X_{\tau_n}=n\;|X_0=u).$$ Hence
\[\frac{1}{M^{\rho-1}}P(X_{\tau_n}< M\;|X_0=u) + \frac{1}{n^{\rho-1}}P(X_{\tau_n}=n\;|X_0=u)\le  \frac{1}{u^{\rho-1}}\;.\]
Therefore
\[P(X_{\tau_n}< M\;|X_0=u)\le \left(\frac{M}{u}\right)^{\rho-1}.\]
Let $n$ go to infinity, we have
$$\psi_M(u)\le\left(\frac{M}{u}\right)^{\rho-1}.$$
Since $\psi(u)\le \psi_M(u)$, we have $$\psi(u)\le\left(\frac{M}{u}\right)^{\rho-1} \;\; \forall \;u\ge M.$$

The cases for $\rho<1$ and $\rho=1$ follow from the next two lemmas.
\begin{lemma}\label{A}
Consider the model given by (\ref{C_L}) and assume that  $\sigma> 0$ and  $\rho<1$. Then there exists $u^{*}>M,$ such that  $$\psi_{u^*}(u)=1,\;\; \forall\; u\ge u^{*}.$$
\end{lemma}
 Proof. Let $F(x)=x^{\alpha}, x>M,$ where $ 0<\alpha<1-\rho .$ Applying It\^{o}'s formula, we have
\begin{align*}
F(X_t)-F(X_0)&=\int_0^t\alpha(X_s)^{\alpha-1}(aX_s+c)\;ds+
\int_0^t\alpha(X_s)^{\alpha-1}\sigma X_s dW_s\\
&+\frac{1}{2}\int_0^t\alpha(\alpha-1)(X_s)^{\alpha-2}\sigma^2{X_s}^2\;ds\\
&+
\int_0^{t^+}\int_0^{M}(X_{s^-}-x)^{\alpha}-(X_{s^-})^{\alpha}\;N_p(dsdx).
\end{align*}
Hence
\begin{align*}
F(X_t)&=F(X_0)+\mbox{ mart. }+\int_0^t\alpha(X_s)^{\alpha-1}(aX_s+c)\;ds\\
&+\frac{1}{2}\int_0^t\alpha(\alpha-1)(X_s)^{\alpha-2}\sigma^2{X_s}^2\;ds\\
&+\int_0^{t^+}\int_0^{M}(X_{s^-}-x)^{\alpha}-(X_{s^-})^{\alpha}\;\tilde{N}_p(dsdx)\\
&\le F(X_0)+\mbox{ mart. }+\alpha\int_0^t (X_s)^{\alpha}\left(\frac{\sigma^2}{2}(\rho+\alpha-1)+cX_s^{-1}\right) ds,
\end{align*}
$\forall\; t\ge 0.$ The above inequality holds because $(X_{s^-}-x)^{\alpha}\le (X_{s^-})^{\alpha},\; \forall X_{s^-}\ge M.$

 Let $u^*=\max(M, 2c/{\sigma^2(1-\rho-\alpha)})$. We consider the process $X_t$ on $[u^*, n)$, where $n$ is an integer ($>u^*$), and let $$\tau_n=\inf\{t>0:\; X_t \not \in [u^*, n)\}$$ be the first exit time from the interval $[u^*, n)$. Then
\begin{eqnarray}\label{Expect}
F(X_{\tau_n})\le F(X_0)+\mbox{ mart. }
\end {eqnarray}
Taking expectation on both sides of the above inequality, and by the Optional Stopping Theorem, we have
\[E[X_{\tau_n}^{\alpha}]\le u^\alpha.\]
Since $F(x)$ is increasing, we have
$$E[F(X_{\tau_n})]\ge(u^*-M)^{\alpha}P(X_{\tau_n}< u^*\;|X_0=u) + n^{\alpha}P(X_{\tau_n}=n\;|X_0=u) .$$
Hence
\[(u^*-M)^{\alpha}P(X_{\tau_n}< u^*\;|X_0=u) + n^{\alpha}P(X_{\tau_n}=n\;|X_0=u)\le u^{\alpha}.\]
 Therefore
\[P(X_{\tau_n}=n\;|X_0=u)\le \left(\frac{u}{n}\right)^{\alpha}.\]
Let $n$ go to infinity, we have
$$\psi_{u^*}(u)=1-\lim_{n\rightarrow \infty} P(X_{\tau_n}=n\;|X_0=u)\ge 1-\lim_{n\rightarrow \infty}\left(\frac{u}{n}\right)^{\alpha}=1,\;\;\forall\; u\ge u^*.$$

\begin{lemma}\label{B}
Consider the model given by (\ref{C_L}) and assume that  $\sigma> 0$ and  $\rho=1$. Then there exists $u^{*}>M+3,$ such that  $$\psi_{u^*}(u)=1\;\; \forall\; u\ge u^{*}.$$
\end{lemma}
 Proof. Let $F(x)=\ln\ln x,\; x>M.$  Applying It\^{o}'s formula, we have
\begin{align*}
F(X_t)-F(X_0)&=\int_0^t(X_s\ln X_s)^{-1}(aX_s+c)\;ds+\int_0^t(X_s\ln X_s)^{-1}\sigma X_s dW_s\\
&+\frac{1}{2}\int_0^t(-\ln{X_s}-1)(X_s\ln X_s)^{-2}\sigma^2{X_s}^2\;ds\\
&+\int_0^{t^+}\int_0^{M}[\ln\ln(X_{s^-}-x)-\ln\ln X_{s^-}]\;N_p(dsdx).
\end{align*}
Hence
\begin{align*}
F(X_t)&=F(X_0)+\mbox{ mart. }+\int_0^t(X_s\ln X_s)^{-1}(aX_s+c)\;ds\\
&+\frac{1}{2}\int_0^t(-\ln X_s-1)(X_s\ln X_s)^{-2}\sigma^2{X_s}^2\;ds\\
&+\int_0^{t^+}\int_0^{M}[\ln\ln(X_{s^-}-x)-\ln\ln X_{s^-}]\;\tilde{N}_p(dsdx)\\
&\le F(X_0)+\mbox{ mart. }+\int_0^t \left(cX_s^{-1}-\frac{\sigma^2}{2\ln X_s}\right)(\ln X_s)^{-1} ds.
\end{align*}
 The above inequality holds because $\ln\ln(X_{s^-}-x)\le \ln\ln X_{s^-},\; \forall X_{s^-}\ge M.$

 Now let $\tilde{u}$ be the solution of $\sigma^2 x=2c\ln x$, and  $u^*=\max(M+3,\tilde{u})$. We consider the process $X_t$ on $[u^*, n)$, where $n$ is an integer ($>u^*$), and let $$\tau_n=\inf\{t>0:\; X_t \not \in [u^*, n)\}$$  be the first exit time from the interval $[u^*, n)$. Then we have
\begin{eqnarray}
F(X_{\tau_n})\le F(X_0)+\mbox{ mart. }
\end{eqnarray}
Taking expectation on both sides of the above inequality, and by the Optional Stopping Theorem, we have
\[E[\ln\ln X_{\tau_n}]\le \ln\ln u.\]
Since $F(x)$ is increasing, we have
\begin{align*}
E[\ln\ln X_{\tau_n}] &\ge \ln\ln (u^*-M)P(X_{\tau_n}< u^*-M\;|X_0=u)\\
&+ \ln\ln n P(X_{\tau_n}=n\;|X_0=u). \nonumber
\end{align*}
Hence
\[\ln\ln (u^*-M)P(X_{\tau_n}< u^*-M\;|X_0=u) + \ln\ln n P(X_{\tau_n}=n\;|X_0=u)\le \ln\ln u.\]
 Therefore
\[P(X_{\tau_n}=n\;|X_0=u)\le \frac{\ln\ln u}{\ln\ln n}\;.\]
Let $n$ go to infinity, we have
$$\psi_{u^*}(u)=1-\lim_{n\rightarrow \infty} P(X_{\tau_n}=n\;|X_0=u)\ge 1-\lim_{n\rightarrow \infty}\frac{\ln\ln u}{\ln\ln n}=1,\;\;\forall \; u\ge u^*.$$

\section{Constantinescu and Thommann's Conjecture}
 In this section, we will prove  that the Constantinescu and Thommann's Conjecture is true.

\begin{lemma}\label{yang}
Let $u^*>0$ be any positive real number. Let $M < \infty$ be an essential range for $\xi_1$.  Suppose $\psi_{u^*}(u)=1$, for all $u\ge u^*.$ Then
  $$\psi_K (u)=1,\;\; \forall\;  u\ge K=\max(u^*-\frac{M}{2}, 0).$$
\end{lemma}

\begin{remark}
 $u^*>0$ in the above Lemma is any positive real number, it needs not be the one defined in Lemma \ref{A} or Lemma \ref{B}.
\end{remark}

Proof. Our first step is to show that for any $0<C_1<1$, there exists a $ \beta_0=\beta_0(M, C_1)$ such that $P\left(X_t \le u^*+\frac{M}{8},\; \forall\;0\le t\le \beta_0\;|\;X_0=u\right)\ge C_1>0,$ for all $u^*\ge u\ge K.$

 Let $Y_t, V_t$ be the same as in lemma \ref{Decreasing}, and $X_t=Y_t^{-1}{V_t}^u$ the solution of (\ref{C_L1}).  Define ${Z_t}^{u^*}=Y_t^{-1}\left(u^*+c\int_0^tY_s\;ds\right)$. Since $d{Z_t}^{u^*}=(aX_t+c)dt+\sigma X_t dW_t$, ${Z_t}^{u^*}$ is a diffusion process. By the continuity of ${Z_t}^{u^*}$, $\forall \; \epsilon>0,\; $ we have
$$P\left(\lim_{\beta\rightarrow 0}\; \sup_{{0\le s\le \beta}}|{Z_s}^{u^*}-u^*|<\epsilon \;\right)=1.$$
Hence for the same $\epsilon>0$ and $\;\forall\; 0<C_1<1.$ $\exists \;\beta_0=\beta_0(\epsilon, C_1)>0,\; s.t.$ \; $$P\left(\sup_{{0\le s\le \beta_0}}|{Z_s}^{u^*}-u^*|<\epsilon\;\right)\ge C_1>0,$$

 In particular, choose $\epsilon=\frac{M}{8},\; \exists \; \beta_0=\beta_0(M, C_1)>0, s.t.$
$$P \left({Z_t}^{u^*} \le u^*+\frac{M}{8},\; \forall\;0\le t\le \beta_0\;\right)\ge C_1>0.$$
 Define ${Z_t}^u=Y_t^{-1}\left(u+c\int_0^tY_s\;ds\right)$, then ${Z_t}^{u^*}\ge{Z_t}^u\ge X_t, \; \forall \; t\ge0$, and
\begin{align*}
P\left(X_t \le u^*+\frac{M}{8},\; \forall\;0\le t\le \beta_0\;|\;X_0=u\right)&\ge P\left({Z_t}^u \le u^*+\frac{M}{8},\; \forall\;0\le t\le \beta_0\;\right)\\
&\ge P \left({Z_t}^{u^*} \le u^*+\frac{M}{8},\; \forall\;0\le t\le \beta_0\;\right)\\
&\ge C_1>0,
\end{align*}
 $\; \forall \;K\le u\le u^*.$\smallskip

 Let $\delta$ be the time that the first jump occurs. Our next step is to show that there exists   $C_2=C_2(C_1,M)>0$ such that $$ P \left(X_{\delta} <K\;|\;X_0=u\right)\ge C_2>0,\; \forall \;K\le u\le u^*.$$

 Notes that  $\; \forall \;K\le u\le u^*,$
\begin{align*}
&P \left(X_t \le u^*+\frac{M}{8},\; \forall\;0\le t\le \beta_0,\;\delta<\beta_0, \;\xi_1>\frac{3M}{4}\;|\;X_0=u\right)\\
&= P \left(X_t \le u^*+\frac{M}{8},\; \forall\;0\le t\le \beta_0\;|\;X_0=u\right)P \left(\delta<\beta_0\right)P \left(\xi_1>\frac{3M}{4}\right) \\
&\ge C_1P \left(\delta<\beta_0\right)P \left(\xi_1>\frac{3M}{4}\right)=C_2>0,
\end{align*}
since $M$ is an essential range of $\xi_1$ and therefore $P( \xi_1 > \frac{3M}{4}) >0$.
  On the other hand,
\begin{align*}
&P \left(X_t \le u^*+\frac{M}{8},\; \forall\;0\le t\le \beta_0,\;\delta<\beta_0, \;\xi_1>\frac{3M}{4}\;|\;X_0=u\right)\\
&\le P\left(X_t \le u^*+\frac{M}{8},\; \forall\;0\le t< \delta,\;\delta<\beta_0,\;\xi_1>\frac{3M}{4}\;|\;X_0=u\right) \\
&\le P \left(X_{\delta} \le u^*+\frac{M}{8}-\frac{3M}{4}=u^*-\frac{5M}{8}<u^*-\frac{M}{2}\le K\;|\;X_0=u\right).
\end{align*}
 Hence $$ P \left(X_{\delta} <K\;|\;X_0=u\right)\ge C_2>0,\; \forall \;K\le u\le u^*.$$

Our final step is to show that $$\psi_K (u)=1,\;\; \forall\;  u\ge K=\max(u^*-\frac{M}{2}, 0).$$
Define \[T_1=\left\{
  \begin{array}{ll}
    \inf\{t>\delta,\; X_t\le u^*\}, & if\; X_{\delta}\ge K \\
    \quad  & \quad \\
    \infty, & if\;
 X_{\delta}< K.
  \end{array}
\right.
\]

Note that the infimum of an empty set is $\infty$. But by the assumption $\psi_{u^*}(u)=1$, for all $u\ge u^*,$ we have $T_1=\infty$ if and only if $X_{\delta}< K$.
 Let $E=\{X_t\ge K,\; \forall\; 0\le t<\infty\}$ and $\theta_s$ be the shift operator, then
\begin{align*}
P(E|\;X_0=u^*)&=E[1_E1_{T_1<\infty}\;|\;X_0=u^*]+E[1_E1_{T_1=\infty}\;|\;X_0=u^*]\\
&=E[1_E1_{T_1<\infty}\;|\;X_0=u^*]\\
&=E[1_{T_1<\infty} \theta_{T_1}[1_E]\;|\;X_0=u^*].
\end{align*}
In what follows, we denote $E_x[1_E]=E[1_E|\;X_0=x]$. By the strong Markov property of $X_t$, we have
\begin{align*}
E[1_{T_1<\infty} \theta_{T_1}[1_E]\;|\;X_0=u^*]&=E[1_{T_1<\infty} E_{X_{T_1}}[1_E]\;|\;X_0=u^*]\\
&\le E\left[1_{T_1<\infty} E_{u^*}[1_E]\;|\;X_0=u^*\right]\\
&= E[1_{T_1<\infty}\;|\;X_0=u^*]E_{u^*}[1_E]\\
&\le (1-C_2)E[1_E\;|\;X_0=u^*]\\
&=P(E|\;X_0=u^*)(1-C_2).
\end{align*}

\noindent The first inequality holds since $K\le X_{T_1}\le u^*$ on $\{T_1<\infty\}$. Hence we have
$$P(E|\;X_0=u^*)=P(E|\;X_0=u^*)(1-C_2).$$
Therefore
$P(E|\;X_0=u^*)=0$, i.e. $\psi_K(u^*)=1$. Since $u\le u^*$, by Lemma \ref{Decreasing}, $$\psi_{K}(u)\ge \psi_K(u^*)=1.$$
The proof is completed.
\begin{theorem}\label{Conj}
Consider the model given by (\ref{C_L}) and assume that  $\sigma> 0$ and  $\rho\le 1$. Suppose also the jump distribution is bounded by an essential range $M>0$. Then $$\psi(u)=1,\;\; \forall \;u\ge 0.$$
\end{theorem}

 Proof. By Lemma \ref{A}, Lemma \ref{B} and the strong Markov property of $X_t$, it is sufficient to show that $$\psi(u)=1,\;\; \forall\; 0\le u\le u^*.$$

 By Lemma \ref{yang},
  $\psi_{K_1} (u)=1,\;\; \forall\;  u\ge K_1=\max(u^*-\frac{M}{2}, 0).$
   Applying Lemma \ref{yang} again, with $u^*$ replaced by $K_1$, we have  $$\psi_{K_2} (u)=1,\;\; \forall\;  u\ge K_2=\max(K_1-M, 0)=\max(u^*-2\frac{M}{2}, 0).$$

 Repeating this argument $N=\lceil\frac{2u^*}{M}\rceil$ times, we have
 $$\psi_{K_N} (u)=1,\;\; \forall\;  u\ge K_N=\max(u^*-N\frac{M}{2}, 0)=0,$$
i.e.,
$$\psi(u)=1,\;\; \forall \;u\ge 0.$$

 Next, we will prove the conjecture true without assuming a cap on the claim size.
\begin{theorem}
Consider the model given by (\ref{C_L}) and assume that  $\sigma> 0$ and  $\rho\le 1$. Then $$\psi(u)=1,\;\; \forall \;u\ge 0.$$
\end{theorem}
Proof. Let $M>0$ be a large constant, define \[\hat{\xi_i}=\left\{
  \begin{array}{ll}
    \xi_i, & if\; \xi_{i}\le M \\
    \quad  & \quad \\
    M, & if\;
 \xi_{i}>M,
 \end{array}
\right.
\]
and $\hat{P_t}=\sum_{j=1}^{N(t)} \hat{\xi_j}\;.$ Let $Y_t, V_t$ be the same as in Lemma \ref{Decreasing}, and $X_t=Y_t^{-1}{V_t}^u$ be the solution of (\ref{C_L1}). Define $$Z_t=Y_t^{-1}\left(u+c\int_0^tY_s\;ds-\int_0^tY_s\;d\hat{P_s}\right),$$ then $Z_t\ge X_t, \; \forall \; t\ge0$. Hence
\begin{eqnarray}
\psi(u)=P(X_t<0, \mbox{ for some } \; 0< t < \infty\;|\;X_0=u)\\
\ge P(Z_t<0, \mbox{ for some } \; 0< t < \infty\;|\;Z_0=u).
\end{eqnarray}
On the other hand, since $dZ_t=(aX_t+c)dt+\sigma X_t dW_t-d\hat{P_t}$,  $Z_t$ satisfies (\ref{C_L}) with bounded claim size distribution. Hence, by Theorem \ref{Conj},
$$P(Z_t<0, \mbox{ for some } \; 0< t < \infty\;|\;Z_0=u)=1, \; \forall \;u\ge 0.$$ Therefore
$$\psi(u)=1,\;\; \forall \;u\ge 0.$$

\end{document}